\documentclass{amsproc}

\def \Z{\mathbb Z}
\def \C{\mathbb C}
\def \R{\mathbb R}

\def \N{\mathbb N}

\def \End{{\rm End}\;}

\def \Aut{{\rm Aut}}
\def \Hom{{\rm Hom}}

\def \<{\langle} 
\def \>{\rangle}

\def \be{\begin{equation}\label}
\def \ee{\end{equation}}
\def \bex{\begin{example}\label}
\def \eex{\end{example}}
\def \bl{\begin{lem}\label}
\def \el{\end{lem}}
\def \bt{\begin{thm}\label}
\def \et{\end{thm}}
\def \bp{\begin{prop}\label}
\def \ep{\end{prop}}
\def \br{\begin{rem}\label}
\def \er{\end{rem}}
\def \bc{\begin{coro}\label}
\def \ec{\end{coro}}
\def \bd{\begin{de}\label}
\def \ed{\end{de}}

\newtheorem{thm}{Theorem}[section]
\newtheorem{prop}[thm]{Proposition}
\newtheorem{coro}[thm]{Corollary}

\newtheorem{example}[thm]{Example}
\newtheorem{lem}[thm]{Lemma}
\newtheorem{rem}[thm]{Remark}
\newtheorem{de}[thm]{Definition}

\numberwithin{equation}{section}

\begin{document}

\title{Twisted modules 
and quasi-modules for vertex operator algebras}

\author{Haisheng Li}
\address{Department of Mathematics, Harbin Normal University, Harbin,
China}
\curraddr{Department of Mathematical Sciences,
Rutgers University, Camden, NJ 08102}
\email{hli@camden.rutgers.edu}
\thanks{The author was supported in part by an NSA Grant.}


\dedicatory{To James Lepowsky and Robert Wilson with
admiration and appreciation.}


\begin{abstract}
We use a result of Barron, Dong and Mason 
to give a natural isomorphism between the category of
twisted modules and the category of quasi-modules of a certain type
for a general vertex operator algebra.
\end{abstract}

\maketitle

\section{Introduction}
In the theory of vertex operator algebras, for a vertex
operator algebra $V$, in addition to the notion of $V$-module 
one has the notion
of $\sigma$-twisted $V$-module where $\sigma$ is a finite order
automorphism of $V$.  For a $V$-module $W$, each element
$v\in V$ is represented by a vertex operator
$$Y_{W}(v,x)\in \Hom (W,W((x)))\subset (\End W)[[x,x^{-1}]],$$
where these vertex operators are mutually local in the sense that 
for $u,v\in V$, there
exists a nonnegative integer $k$ such that
$$(x_{1}-x_{2})^{k}[Y_{W}(u,x_{1}),Y_{W}(v,x_{2})]=0.$$
Twisted modules were first introduced and used by 
Frenkel, Lepowsky and Meurman in their construction of 
the moonshine module vertex operator algebra $V^{\natural}$ 
(see \cite{lep1}, \cite{flm}). 
Let $V$ be a vertex operator algebra and let $\sigma$ be an
automorphism of order $N$. 
For a $\sigma$-twisted $V$-module $W$
(\cite{lep1}, \cite{flm}, \cite{ffr}, \cite{dong}),
each element $v$ of $V$ is represented by a twisted vertex operator
$$Y_{W}(v,x)\in \Hom (W,W((x^{1/N})))
\subset (\End W)[[x^{1/N},x^{-1/N}]],$$
where these twisted vertex operators are also mutually local.

In a recent work \cite{li-ga}, to associate certain
(infinite-dimensional) Lie algebras with vertex algebras, we studied
what we called quasi local vertex operators (cf. \cite{gkk}).  Let $W$ be
any vector space. A subset $S$ of $\Hom (W,W((x)))$ is said to be
quasi local if for any $a(x), b(x)\in S$ there exists a nonzero
polynomial $p(x_{1},x_{2})$ such that
$$p(x_{1},x_{2})a(x_{1})b(x_{2})=p(x_{1},x_{2})b(x_{2})a(x_{1}).$$
It was proved therein that any maximal quasi local subspace 
has a natural vertex algebra structure and any quasi local subset 
generates a vertex algebra.
This particular result generalizes the main result of
\cite{li-local}, which states that for any vector space $W$,
any set of mutually local vertex operators on $W$ generates 
a vertex algebra with $W$ as a natural module.
However, the space $W$ under the natural action is not a module 
for vertex algebras generated by quasi local vertex operators on $W$,
though a certain weaker version of Jacobi identity was proved to hold.
This motivated us to introduce a new notion of quasi module
for a vertex algebra. For a quasi module $W$ for a vertex algebra $V$, 
each element $v$ of $V$ is represented by a vertex operator 
$$Y_{W}(v,x)\in \Hom (W,W((x)))\subset (\End W)[[x,x^{-1}]]$$ 
and the vertex operators $Y_{W}(v,x)$ for $v\in V$ form a 
quasi local subspace. 

On twisted modules for vertex operator algebras there is a conceptual
work \cite{bdm}, in which for any vertex operator algebra $V$ and for
any positive integer $k$, a canonical isomorphism was established
between the category of $V$-modules and the category of twisted
modules for the tensor product vertex operator algebra $V^{\otimes k}$
with respect to permutation automorphisms.  In \cite{bdm}, a central
role was played by the geometric change-of-coordinate $x=z^{k}$. It
has been known (\cite{zhu}, [H1-3], cf. \cite{lep2}) that for any
vertex operator algebra $V$ and for any $f(z)\in z\C[[z]]$ with
$f'(0)\ne 0$, the change-of-coordinate $x=f(z)$ gives rise to a
``new'' vertex operator algebra structure on $V$, which was proved to
be isomorphic to $V$.  A special change-of-coordinate played a very
important role in the study of modular invariance of graded characters
(\cite{zhu}, \cite{dlm-modular}).

It has been well known (cf. \cite{fz}) that (untwisted) affine Lie
algebras together with their highest weight modules can be naturally
associated with vertex operator algebras and their modules.
Furthermore, twisted affine Lie algebras together with their highest
weight modules (see \cite{kac}) can be associated with twisted modules
for those vertex operator algebras (cf. \cite{flm},
\cite{li-twisted}).  On the other hand, it was proved in \cite{li-ga}
that twisted affine Lie algebras, which are represented in a different
form, together with their highest weight modules, can be naturally
associated with quasi-modules for the vertex operator algebras
associated with the untwisted affine Lie algebras.  This suggests that
there exist a natural connection between twisted modules and quasi
modules for a general vertex operator algebra.

The main purpose of this paper is to give a natural connection between
twisted modules and quasi-modules for a general vertex operator
algebra.  Indeed, the goal has been achieved by using \cite{bdm},
thanks to Barron, Dong and Mason for their beautiful work.  
What we have proved is that the same change-of-coordinates, used by
Barron, Dong and Mason, give rise to a natural isomorphism between the
category of twisted modules and the category of quasi-modules of a
certain special type.

We thank Yi-Zhi and Kailash for organizing this great conference, in
honor of Professors James Lepowsky and Robert Wilson.  I am very grateful for
having Jim and Robert as teachers and as friends as well.

\section{Twisted modules and quasi-modules}
We here present the main result, a natural isomorphism
between the category of twisted modules 
and the category of quasi-modules of a certain type
for a general vertex operator algebra.

First, we recall the definitions of the notions of twisted module and
quasi-module.  Let $V=\coprod_{n\in \Z}V_{(n)}$ be a vertex operator
algebra, fixed throughout this section. For the definition and basic
properties we refer to \cite{flm} and \cite{fhl}.  Let $\sigma$ be an
automorphism of $V$ of order $N$ (a positive integer).  Then
$V=\coprod_{j=0}^{N-1}V^{j}$, where $V^{j}=\{ u\in V\;|\; \sigma
(u)=\omega_{N}^{j}u\}$ and $\omega_{N}=\exp (2\pi \sqrt{-1}/N)$, the
principal primitive $N$-th root of unity.

A {\em weak $\sigma$-twisted $V$-module} 
(cf. \cite{lep1}, \cite{flm}, \cite{ffr}, \cite{dong}, \cite{dlm-twisted})
is a vector space $W$ equipped with a linear map
\begin{eqnarray}
Y_{W}:& & V\rightarrow \Hom (W,W((x^{\frac{1}{N}})))\subset
(\End W)[[x^{\frac{1}{N}},x^{-\frac{1}{N}}]]\nonumber\\
& &v\rightarrow Y_{W}(v,x)
\end{eqnarray}
such that 
$$Y_{W}({\bf 1},x)=1_{W}\;\mbox{(the identity operator on $W$)},$$
and  for $u, v\in V$,
\begin{eqnarray}\label{etwistedmodule-def}
& &x_{0}^{-1}\delta\left(\frac{x_{1}-x_{2}}{x_{0}}\right)
Y_{W}(u,x_{1})Y_{W}(v,x_{2})
-x_{0}^{-1}\delta\left(\frac{x_{2}-x_{1}}{-x_{0}}\right)
Y_{W}(v,x_{2})Y_{W}(u,x_{1})\nonumber\\
& &\ \ \ \ \ \ \ =\frac{1}{N}\sum_{r=0}^{N-1}x_{1}^{-1}\delta\left(
\omega_{N}^{r}\left(\frac{x_{2}+x_{0}}{x_{1}}\right)^{\frac{1}{N}}\right)
Y_{W}(Y(\sigma^{r}u,x_{0})v,x_{2})
\end{eqnarray}
(the {\em twisted Jacobi identity}).
Note that as a convention, 
for $\alpha\in \R$, the expressions $(x_{1}\pm x_{2})^{\alpha}$ 
are understood as the formal series 
in the nonnegative integral powers of the second variable $x_{2}$. That is,
$$(x_{1}\pm x_{2})^{\alpha}=\sum_{i\in \N}\binom{\alpha}{i}(\pm 1)^{i}
x_{1}^{\alpha-i}x_{2}^{i}
\in x_{1}^{\alpha}\R[x_{1}^{-1}][[x_{2}]].$$

If $u\in V^{j}$ with $0\le j\le N-1$, 
the twisted Jacobi identity becomes
\begin{eqnarray}
& &x_{0}^{-1}\delta\left(\frac{x_{1}-x_{2}}{x_{0}}\right)
Y_{W}(u,x_{1})Y_{W}(v,x_{2})
-x_{0}^{-1}\delta\left(\frac{x_{2}-x_{1}}{-x_{0}}\right)
Y_{W}(v,x_{2})Y_{W}(u,x_{1})\nonumber\\
& &\ \ \ \ \ \ \ =x_{1}^{-1}\delta\left(\frac{x_{2}+x_{0}}{x_{1}}\right)
\left(\frac{x_{2}+x_{0}}{x_{1}}\right)^{-\frac{j}{N}}Y_{W}(Y(u,x_{0})v,x_{2}).
\end{eqnarray}
By taking $v={\bf 1}$, one obtains
\begin{eqnarray}
Y_{W}(u,x)\in x^{\frac{j}{N}}\Hom (W,W((x))).
\end{eqnarray}
This particular property amounts to
\begin{eqnarray}
Y_{W}(\sigma u,x)=\lim_{x^{1/N}\rightarrow \omega_{N}x^{1/N}}Y_{W}(u,x).
\end{eqnarray}

\br{rinverse}
{\em Note that the above defined notion of $\sigma$-twisted
  $V$-module, which is the one defined in 
\cite{dlm-modular} and \cite{bdm}, 
corresponds to the notion of $\sigma^{-1}$-twisted $V$-module
in \cite{dlm-twisted}.}
\er

The twisted Jacobi identity is equivalent to 
the following {\em weak commutativity and associativity} 
(\cite{dl}, \cite{li-twisted}):
For $u,v\in V$, there exists a nonnegative integer $k$ such that
\begin{eqnarray}
(x_{1}-x_{2})^{k}[Y_{W}(u,x_{1}),Y_{W}(v,x_{2})]=0,
\end{eqnarray}
and for $u\in V^{j},\; v\in V,\; w\in W$, $0\le j\le N-1$,
there exists a nonnegative integer $l$ such that
\begin{eqnarray}
& &(x_{0}+x_{2})^{l-j/N}Y_{W}(u,x_{0}+x_{2})Y_{W}(v,x_{2})w\\
&=&(x_{2}+x_{0})^{l-j/N}Y_{W}(Y(u,x_{0})v,x_{2})w.\nonumber
\end{eqnarray}

{}From now on we fix an automorphism $\sigma$ of order
$N$ for the fixed vertex operator algebra $V$. Set
\begin{eqnarray}
G=\<\sigma\>\subset \Aut (V)
\mbox{ (the full automorphism group of $V$)}.
\end{eqnarray}
Let $\phi: G\rightarrow \C^{\times}$ be the (injective) 
group homomorphism defined by
$\phi(\sigma^{j})=\omega_{N}^{j}$ for $j=0,\dots,N-1$.

A {\em quasi $V$-module} \cite{li-ga} is a vector space $W$ 
equipped with a linear map
\begin{eqnarray*}
Y_{W}: & &V\rightarrow \Hom (W,W((x)))\subset (\End W)[[x,x^{-1}]],\\
& &v\rightarrow Y_{W}(v,x)
\end{eqnarray*}
such that
$$Y_{W}({\bf 1},x)=1_{W},$$
and  for $u,v\in V$, there exists a nonzero polynomial $p(x_{1},x_{2})$
such that
\begin{eqnarray}\label{equasi-jacobi}
& &x_{0}^{-1}\delta\left(\frac{x_{1}-x_{2}}{x_{0}}\right)
p(x_{1},x_{2})Y_{W}(u,x_{1})Y_{W}(v,x_{2})\nonumber\\
& &\hspace{2cm}-x_{0}^{-1}\delta\left(\frac{x_{2}-x_{1}}{-x_{0}}\right)
p(x_{1},x_{2})Y_{W}(v,x_{2})Y_{W}(u,x_{1})\nonumber\\
&=&x_{2}^{-1}\delta\left(\frac{x_{1}-x_{0}}{x_{2}}\right)p(x_{1},x_{2})
Y_{W}(Y(u,x_{0})v,x_{2})
\end{eqnarray}
(the {\em quasi Jacobi identity}).

\bd{dgphi-module}
{\em A quasi $V$-module $(W,Y_{W})$ is called 
a {\em $(G,\phi)$-quasi $V$-module}  
if for any $u,v\in V$, there exists a nonnegative integer $k$ such that
(\ref{equasi-jacobi}) holds with
$p(x_{1},x_{2})=(x_{1}^{N}-x_{2}^{N})^{k}$
and such that
\begin{eqnarray}
Y_{W}(\phi(g)^{-L(0)}g(u),x)=Y_{W}(u,\phi(g)x)\;\;\;\mbox{ for }g\in G,\; u\in V.
\end{eqnarray}}
\ed

Follow \cite{bdm} to define $a_{n}\in \C$ for $n\in \Z_{+}$ by
\begin{eqnarray}
\exp \left(-\sum_{n\in \Z_{+}}a_{n}x^{n+1}\frac{d}{dx}\right)\cdot x
=\frac{1}{N}(1+x)^{N}-\frac{1}{N},
\end{eqnarray}
where $\Z_{+}$ denotes the set of positive integers, and then set
\begin{eqnarray}
\Delta_{N}(x)=\exp\left(\sum_{n\in \Z_{+}}a_{n}x^{-n/N}L(n)\right)
N^{-L(0)}x^{(1/N-1)L(0)},
\end{eqnarray}
an invertible element of $(\End V)[[x^{1/N},x^{-1/N}]]$. We have
\begin{eqnarray}
\Delta_{N}(x^{N})^{-1}
=(Nx^{N-1})^{L(0)}\exp\left(-\sum_{n\in \Z_{+}}a_{n}x^{-n}L(n)\right).
\end{eqnarray}
As in \cite{bdm} we shall also heavily use the expression 
$\Delta_{N}(x^{N})^{-1}$. For convenience we set
\begin{eqnarray}
\Phi(x)=\Delta_{N}(x^{N})^{-1}\in \Hom (V,V[x,x^{-1}]).
\end{eqnarray}

The following result was proved in \cite{bdm}:

\bp{pbdm}
For any $u\in V$, we have
\begin{eqnarray}\label{ebdm1}
\ \ \ \ \Delta_{N}(x)Y(u,x_{0})\Delta_{N}(x)^{-1}
=Y(\Delta_{N}(x+x_{0})u,(x+x_{0})^{1/N}-x^{1/N})
\end{eqnarray}
in $(\End V)[[x_{0}^{\pm 1},x^{\pm \frac{1}{N}}]]$, and
\begin{eqnarray}\label{ebdm2}
\Phi(x)Y(u,x_{0})
=Y(\Phi(x+x_{0})u,(x+x_{0})^{N}-x^{N})\Phi(x)
\end{eqnarray}
in $(\End V)[[x_{0}^{\pm 1},x^{\pm 1}]]$.
\ep

\br{rsubstitutions1}
{\em 
Recall from \cite{bdm} the interpretation of the formal variable notations
in Proposition \ref{pbdm}.
First, for any nonzero $\alpha\in \frac{1}{N}\Z$, under the convention we have
$$(x+x_{0})^{\alpha}-x^{\alpha}
=\sum_{i\in \Z_{+}}\binom{\alpha}{i}x^{\alpha-i}x_{0}^{i}
=x^{\alpha-1}x_{0}(\alpha +x_{0}f),$$
where 
$f=\sum_{i\ge 2}\binom{\alpha}{i}x^{1-i}x_{0}^{i-2}
\in \R[x^{-1}][[x_{0}]]$.
For $n\in \Z$, it is understood that
\begin{eqnarray*}
\left((x+x_{0})^{\alpha}-x^{\alpha}\right)^{n}
=x^{n(\alpha -1)}x_{0}^{n}
\sum_{i\in \N}\binom{n}{i}\alpha^{n-i}x_{0}^{i}f^{i}
\in x^{n\alpha}x_{0}^{n}\R[x,x^{-1}][[x_{0}]].
\end{eqnarray*}
Then for $u,v\in V$, we have
\begin{eqnarray*}
Y(u,(x+x_{0})^{\alpha}-x^{\alpha})v
=\sum_{m\in \Z}u_{m}v\left((x+x_{0})^{\alpha}-x^{\alpha}\right)^{-m-1}
\in V[x^{\pm 1/N}]((x_{0})).
\end{eqnarray*}
This explains the formal variable notations in Proposition \ref{pbdm}.
As we shall mention in the following Remark, we shall also use
another (different) substitution. For this purpose, we also write
$$Y(u,z)|_{z=(x+x_{0})^{\alpha}-x^{\alpha},x>>x_{0}}$$
for this particularly defined expression $Y(u,(x+x_{0})^{\alpha}-x^{\alpha})$.
It was showed in \cite{bdm}, page 363, that 
for $h,\alpha\in \frac{1}{N}\Z$,
\begin{eqnarray}\label{ebdm363}
(x_{2}^{\alpha}+z_{0})^{h}|_{z_{0}=(x_{2}+x_{0})^{\alpha}-x_{2}^{\alpha},
x_{2}>>x_{0}}
=(x_{2}+x_{0})^{\alpha h}.
\end{eqnarray}}
\er

{\em Warning:} The following expansion
\begin{eqnarray*}
(z-x^{\alpha})^{n}|_{z=(x+x_{0})^{\alpha}}&=&
\sum_{i\in \N}(-1)^{i}\binom{n}{i}(x+x_{0})^{\alpha (n-i)}x^{\alpha i}\\
&=&\sum_{i\in \N}\sum_{j\in \N}(-1)^{i}
\binom{n}{i}\binom{(n-i)\alpha}{j}x^{n\alpha-j}x_{0}^{j}
\end{eqnarray*}
is an infinite divergent sum if $n<0$.

\br{rsubstitutions}
{\em We shall need a different substitution 
$z=(x+x_{0})^{N}-x^{N}$ for rational powers $z^{\alpha},$ 
$\alpha\in \frac{1}{N}\Z$. Let $p(x_{0},x)=x_{0}^{k}+xq(x_{0},x)$, where
$k$ is a positive integer and $q(x_{0},x)$ is a polynomial. 
We consider the following expansion
$$p(x_{0},x)^{\alpha}=(x_{0}^{k}+xq(x_{0},x))^{\alpha}=\sum_{i\in
    \N}\binom{\alpha}{i}x_{0}^{k(\alpha-i)}x^{i}q(x_{0},x)^{i}\in 
\R[x_{0},x_{0}^{-1}][[x]],$$
and we use the notation
$$z^{\alpha}|_{z=p(x_{0},x),x_{0}>>x}$$
for this particular expansion. That is,
$$z^{\alpha}|_{z=p(x_{0},x),x_{0}>>x}=(x_{0}^{k}+y)^{\alpha}|_{y=xq(x_{0},x)}$$
under the usual convention.
In particular, for $p(x_{0},x)=x_{0}+x$ we have
$$z^{\alpha}|_{z=x_{0}+x,x_{0}>>x}=(x_{0}+x)^{\alpha}.$$ 
For any $h(x_{0},x)\in \R[x_{0},x]$ we  have
\begin{eqnarray}
& &z^{\alpha}|_{z=p(x_{0},x)+xh(x_{0},x),x_{0}>>x}\nonumber\\
&=&\sum_{i\in \N}\binom{\alpha}{i}
x_{0}^{k(\alpha -i)}(xq(x_{0},x)+xh(x_{0},x))^{i}\nonumber\\
&=&\sum_{i,j\in \N}\binom{\alpha}{i}\binom{i}{j}
x_{0}^{k(\alpha -i)}(xq(x_{0}+x))^{i-j}(xh(x_{0},x))^{j}\nonumber\\
&=&\sum_{r,j\in \N}\binom{\alpha}{r+j}\binom{r+j}{j}
x_{0}^{k(\alpha -r-j)}(xq(x_{0}+x))^{r}(xh(x_{0},x))^{j}\nonumber\\
&=&\sum_{r\in\N}\sum_{j\in \N}\binom{\alpha}{j}\binom{\alpha-j}{r}
x_{0}^{k(\alpha -r-j)}(xq(x_{0}+x))^{r}(xh(x_{0},x))^{j}
\nonumber\\
&=&\sum_{j\in \N}\binom{\alpha}{j}
z_{0}^{\alpha-j}|_{z_{0}=p(x_{0},x),x_{0}>>x}(xh(x_{0},x))^{j}\nonumber\\
&=&(z_{0}+xh(x_{0},x))^{\alpha}|_{z_{0}=p(x_{0},x),x_{0}>>x}.
\end{eqnarray}
Using this and (\ref{ebdm363}) we get
\begin{eqnarray}
(z_{0}+x^{N})^{\alpha}|_{z_{0}=(x_{0}+x)^{N}-x^{N},x_{0}>>x}
&=&z^{\alpha}|_{z=(x_{0}+x)^{N},x_{0}>>x}\nonumber\\
&=&(x_{0}^{N}+y)^{\alpha}|_{y=(x_{0}+x)^{N}-x_{0}^{N},x_{0}>>x}\nonumber\\
&=&(x_{0}+x)^{N\alpha}.
\end{eqnarray}}
\er

We shall need the following simple result:

\bl{leasy}
For any $\alpha\in \C^{\times}$,
\begin{eqnarray}
& &\Phi(x)\alpha^{-L(0)}=\alpha^{-NL(0)}\Phi(\alpha x),\\
& &\alpha^{L(0)}\Delta_{N}(x)
=\lim_{x^{1/N}\rightarrow \alpha x^{1/N}}\Delta_{N}(x)\alpha^{NL(0)}
\end{eqnarray}
hold in $\Hom (V,V[x,x^{-1}])$ and 
in $\Hom (V,V[x^{\frac{1}{N}},x^{-\frac{1}{N}}])$, respectively.
\el

\begin{proof} We shall just prove the first identity, as the second will 
follow easily.
Since $[L(0),L(n)]=-nL(n)$ for $n\in \Z$,
it follows that 
$$\alpha^{L(0)}L(n)\alpha^{-L(0)}=\alpha^{-n}L(n).$$
Thus
\begin{eqnarray*}
\alpha^{L(0)}\left(\sum_{n\in \Z_{+}}a_{n}x^{-n}L(n)\right)\alpha^{-L(0)}
&=&\sum_{n\in \Z_{+}}\alpha^{-n}a_{n}x^{-n}L(n)\\
&=&\sum_{n\in \Z_{+}}a_{n}(\alpha x)^{-n}L(n).
\end{eqnarray*}
Using this we obtain
\begin{eqnarray*}
\Phi(x)\alpha^{-L(0)}
&=&(Nx^{N-1})^{L(0)}\alpha^{-L(0)}
\alpha^{L(0)}\exp\left(-\sum_{n\in \Z_{+}}a_{n}x^{-n}L(n)\right)
\alpha^{-L(0)}\\
&=&\alpha^{-NL(0)}(N(\alpha x)^{N-1})^{L(0)}
\exp\left(-\sum_{n\in \Z_{+}}a_{n}(\alpha x)^{-n}L(n)\right)\\
&=&\alpha^{-NL(0)}\Phi(\alpha x),
\end{eqnarray*}
proving the assertion.
\end{proof}

The following is the first half of our main result of this paper:

\bt{tmain1}
Let $(W,Y_{W})$ be a weak $\sigma$-twisted $V$-module. 
For $u\in V$, as in \cite{bdm} set
\begin{eqnarray}
\tilde{Y}_{W}(u,x)=Y_{W}(\Phi(x)u,x^{N})\in (\End W)[[x,x^{-1}]].
\end{eqnarray}
Then
$(W,\tilde{Y}_{W})$ carries the structure of a $(G,\phi)$-quasi $V$-module.
\et

\begin{proof} First, for $u\in V$,  as $\Phi(x)u\in V[x,x^{-1}]$, we have
$$\tilde{Y}_{W}(u,x)\in \Hom (W,W((x))). $$
Second, we have
$$\tilde{Y}_{W}({\bf 1},x)=Y_{W}({\bf 1},x^{N})=1_{W},$$
as $\Phi(x){\bf 1}={\bf 1}$, which is due to the fact that
$L(n){\bf 1}=0$ for $n\ge 0$.

Third, {}from Lemma \ref{leasy}, 
for any $N$-th root of unity $\alpha$, we have
$$\Phi(x)\alpha^{-L(0)}=\Phi(\alpha x).$$
For $g\in G,\; u\in V$, we have
\begin{eqnarray*}
\tilde{Y}_{W}(\phi(g)^{-L(0)}g (u),x)
&=&Y_{W}(\Phi(x)\phi(g)^{-L(0)}g (u),x^{N})\\
&=&Y_{W}(\Phi(\phi(g)x)g (u),x^{N})\\
&=&\lim_{z\rightarrow \phi(g)x}Y_{W}(\Phi(\phi(g)x)u,z^{N})\\
&=&\tilde{Y}_{W}(u,\phi(g)x),
\end{eqnarray*}
noticing that $g\Phi(x)=\Phi(x)g$.

Now it remains to prove the quasi Jacobi identity.
Let $u,v\in V$. Since $\Phi(x)u,\; \Phi(x)v\in V[x,x^{-1}]$,
there exists a nonnegative integer $k$ such that
\begin{eqnarray}\label{e2.26}
z^{k}Y(\Phi(x_{1})u,z)\Phi(x_{2})v \in V[x_{1}^{\pm 1},x_{2}^{\pm 1},z].
\end{eqnarray}
Then
$$(z_{1}-z_{2})^{k}[Y_{W}(\Phi(x_{1})u,z_{1}),
Y_{W}(\Phi(x_{2})v,z_{2})]=0,$$
which yields
\begin{eqnarray}\label{eweak-comm-proof}
(x_{1}^{N}-x_{2}^{N})^{k}[\tilde{Y}_{W}(u,x_{1}),\tilde{Y}_{W}(v,x_{2})]=0.
\end{eqnarray}
Let $r$ be a nonnegative
integer such that  $x^{r}\Phi(x)u\in V[x]$.
Then 
$$(x_{0}+x_{2})^{r}\Phi(x_{0}+x_{2})u
=(x_{2}+x_{0})^{r}\Phi(x_{2}+x_{0})u\in V[x_{0},x_{2}].$$
Furthermore, let $w\in W$. In view of Remark \ref{rsubstitutions}, we have
\begin{eqnarray*}
& &(x_{0}+x_{2})^{r}\tilde{Y}_{W}(u,x_{0}+x_{2})\tilde{Y}_{W}(v,x_{2})w\\
&=&(x_{0}+x_{2})^{r}Y_{W}(\Phi(x_{0}+x_{2})u,(x_{0}+x_{2})^{N})
Y_{W}(\Phi(x_{2})v,x_{2}^{N})w\\
&=&\left((x_{0}+x_{2})^{r}Y_{W}(\Phi(x_{0}+x_{2})u,z_{0}+x_{2}^{N})
Y_{W}(\Phi(x_{2})v,x_{2}^{N})w\right)\\
& &\hspace{6cm}|_{z_{0}=(x_{0}+x_{2})^{N}-x_{2}^{N},x_{0}>>x_{2}}
\end{eqnarray*}
and
\begin{eqnarray*}
& &(x_{2}+x_{0})^{r}\tilde{Y}_{W}(Y(u,x_{0})v,x_{2})w\\
&=&(x_{2}+x_{0})^{r}Y_{W}(\Phi(x_{2})Y(u,x_{0})v,x_{2}^{N})w\\
&=&(x_{2}+x_{0})^{r}
Y_{W}\left(Y(\Phi(x_{2}+x_{0})u,(x_{2}+x_{0})^{N}-x_{2}^{N})
\Phi(x_{2})v,x_{2}^{N}\right)w\\
&=&\left((x_{2}+x_{0})^{r}
Y_{W}\left(Y(\Phi(x_{2}+x_{0})u,z_{0})
\Phi(x_{2})v,x_{2}^{N}\right)w\right)|_{z_{0}
=(x_{2}+x_{0})^{N}-x_{2}^{N},x_{2}>>x_{0}},
\end{eqnarray*}
using (\ref{ebdm2}).

Assume $u\in V^{j}$ with $0\le j\le N-1$. 
There exists a positive integer $l$ such that
\begin{eqnarray*}
& &(z_{0}+x_{2}^{N})^{l-\frac{j}{N}}
Y_{W}((x_{0}+x_{2})^{r}\Phi(x_{0}+x_{2})u,z_{0}+x_{2}^{N})
Y_{W}(\Phi(x_{2})v,x_{2}^{N})w\nonumber\\
&=&(x_{2}^{N}+z_{0})^{l-\frac{j}{N}}
Y_{W}(Y((x_{0}+x_{2})^{r}\Phi(x_{0}+x_{2})u,z_{0})\Phi(x_{2})v,x_{2}^{N})w,
\end{eqnarray*}
which gives
\begin{eqnarray}
& &z_{0}^{k}(z_{0}+x_{2}^{N})^{l-\frac{j}{N}}Y_{W}((x_{0}+x_{2})^{r}\Phi(x_{0}+x_{2})u,z_{0}+x_{2}^{N})
Y_{W}(\Phi(x_{2})v,x_{2}^{N})w\\
&=&z_{0}^{k}(x_{2}^{N}+z_{0})^{l-\frac{j}{N}}
Y_{W}(Y((x_{0}+x_{2})^{r}\Phi(x_{0}+x_{2})u,z_{0})\Phi(x_{2})v,x_{2}^{N})w,
\nonumber
\end{eqnarray}
where $k$ is the nonnegative integer as in (\ref{e2.26}).
Now, we shall perform the substitution
$z_{0}=(x_{2}+x_{0})^{N}-x_{2}^{N},x_{0}>>x_{2}$ 
on both sides.
Notice that the expression on the right-hand side
involves nonnegative integral powers of $z_{0}$ only, so that
the substitutions
$z_{0}=(x_{2}+x_{0})^{N}-x_{2}^{N},x_{0}>>x_{2}$ 
and $z_{0}=(x_{2}+x_{0})^{N}-x_{2}^{N},x_{2}>>x_{0}$ 
agree on the right-hand side.
Performing the substitution
$z_{0}=(x_{2}+x_{0})^{N}-x_{2}^{N},x_{0}>>x_{2}$ 
on both sides and using Remark \ref{rsubstitutions}, 
we obtain
\begin{eqnarray*}
& &(x_{0}+x_{2})^{r+Nl-j}((x_{0}+x_{2})^{N}-x_{2}^{N})^{k}\cdot \\
& &\hspace{2cm}\cdot
Y_{W}(\Phi(x_{0}+x_{2})u,(x_{0}+x_{2})^{N})Y_{W}(\Phi(x_{2})v,x_{2}^{N})w\\
&=&\left((x_{2}+x_{0})^{r}(x_{2}^{N}+z_{0})^{l-\frac{j}{N}}z_{0}^{k}
Y_{W}\left(Y(\Phi(x_{2}+x_{0})u,z_{0})
\Phi(x_{2})v,x_{2}^{N}\right)w\right)\\
& &\hspace{6cm}|_{z_{0}=(x_{2}+x_{0})^{N}-x_{2}^{N},x_{0}>>x_{2}}\\
&=&\left((x_{2}+x_{0})^{r}(x_{2}^{N}+z_{0})^{l-\frac{j}{N}}z_{0}^{k}
Y_{W}\left(Y(\Phi(x_{2}+x_{0})u,z_{0})
\Phi(x_{2})v,x_{2}^{N}\right)w\right)\\
& &\hspace{6cm}|_{z_{0}=(x_{2}+x_{0})^{N}-x_{2}^{N},x_{2}>>x_{0}}\\
&=&(x_{2}+x_{0})^{r+Nl-j}\left((x_{2}+x_{0})^{N}-x_{2}^{N}\right)^{k}
\tilde{Y}_{W}(Y(u,x_{0})v,x_{2})w.
\end{eqnarray*}
Thus
\begin{eqnarray}\label{eweak-assoc-proof}
\lefteqn{(x_{0}+x_{2})^{r+Nl-j}((x_{0}+x_{2})^{N}-x_{2}^{N})^{k}
\tilde{Y}_{W}(u,x_{0}+x_{2})\tilde{Y}_{W}(v,x_{2})w}\\
&=&(x_{0}+x_{2})^{r+Nl-j}((x_{0}+x_{2})^{N}-x_{2}^{N})^{k}
\tilde{Y}_{W}(Y(u,x_{0})v,x_{2})w.\nonumber
\end{eqnarray}
Combining (\ref{eweak-comm-proof}) and (\ref{eweak-assoc-proof})  
we obtain the following quasi Jacobi identity
\begin{eqnarray}
\lefteqn{x_{0}^{-1}\delta\left(\frac{x_{1}-x_{2}}{x_{0}}\right)
(x_{1}^{N}-x_{2}^{N})^{k}\tilde{Y}_{W}(u,x_{1})\tilde{Y}_{W}(v,x_{2})w}
\nonumber\\
& &\hspace{1cm}-x_{0}^{-1}\delta\left(\frac{x_{2}-x_{1}}{-x_{0}}\right)
(x_{1}^{N}-x_{2}^{N})^{k}\tilde{Y}_{W}(v,x_{2})\tilde{Y}_{W}(u,x_{1})w\nonumber\\
&=&x_{2}^{-1}\delta\left(\frac{x_{1}-x_{0}}{x_{2}}\right)(x_{1}^{N}-x_{2}^{N})^{k}
\tilde{Y}_{W}(Y(u,x_{0})v,x_{2})w.
\end{eqnarray}
Therefore, $(W,\tilde{Y}_{W})$ carries the structure of a $(G,\phi)$-quasi $V$-module.
 \end{proof}

\br{rmodule-quasi}
{\em  Let $(W,Y_{W})$ be any weak $V$-module. Then the same proof of
Theorem \ref{tmain1} 
shows that $(W,\tilde{Y}_{W})$ is a $(G,\phi)$-quasi $V$-module.}
\er

Next we present the second half of our main result of this paper.

\bt{tconverse}
Let $(W,Y_{W})$ be a $(G,\phi)$-quasi $V$-module.
For $u\in V$, as in \cite{bdm} set
\begin{eqnarray}
\bar{Y}_{W}(u,x)
=Y_{W}(\Delta_{N}(x)u,x^{\frac{1}{N}})
\in (\End W)[[x^{\frac{1}{N}},x^{-\frac{1}{N}}]].
\end{eqnarray}
Then $(W,\bar{Y}_{W})$ carries the structure of 
a weak $\sigma$-twisted $V$-module.
\et

\begin{proof} For convenience, let us simply use 
$\Delta(x)$ for $\Delta_{N}(x)$ in the proof.
First, for $u\in V$, 
as $\Delta(x)u\in V[x^{1/N},x^{-1/N}]$, we have
$$\bar{Y}_{W}(u,x)=Y_{W}(\Delta(x)u,x^{\frac{1}{N}})
\in \Hom (W,W((x^{\frac{1}{N}}))).$$
Second, 
$$\bar{Y}_{W}({\bf 1},x)=Y_{W}(\Delta(x){\bf 1},x^{\frac{1}{N}})
=Y_{W}({\bf 1},x^{\frac{1}{N}})=1_{W}.$$
Third, for $u\in V$, as $\phi(\sigma)=\omega_{N}$, we have
$Y_{W}(\omega_{N}^{-L(0)}\sigma(u),x)=Y_{W}(u,\omega_{N}x)$.
Using Lemma \ref{leasy}, we get
\begin{eqnarray}\label{elimit-property}
\bar{Y}_{W}(\sigma(u),x)
&=&Y_{W}(\Delta(x)\sigma (u),x^{\frac{1}{N}})\nonumber\\
&=&Y_{W}(\omega_{N}^{-L(0)}\sigma \omega_{N}^{L(0)}
\Delta(x)u,x^{\frac{1}{N}})\nonumber\\
&=&\lim_{z^{\frac{1}{N}}\rightarrow \omega_{N}x^{\frac{1}{N}}}
Y_{W}(\omega_{N}^{L(0)}\Delta(x)u,z^{\frac{1}{N}})\nonumber\\
&=&\lim_{x^{\frac{1}{N}}\rightarrow \omega_{N}x^{\frac{1}{N}}}
Y_{W}(\Delta(x)u,x^{\frac{1}{N}})\nonumber\\
&=&\lim_{x^{\frac{1}{N}}\rightarrow \omega_{N}x^{\frac{1}{N}}}\bar{Y}_{W}(u,x).
\end{eqnarray}
Next, we prove the weak commutativity and the weak associativity, which
amount to the twisted Jacobi identity.
Let $u,v\in V$. As $\Delta(x)u,\; \Delta(x)v
\in V[x^{1/N},x^{-1/N}]$, from the quasi Jacobi identity
there exists a nonnegative integer $k$ such that
\begin{eqnarray}\label{ejacobi-deform}
& &z_{0}^{-1}\delta\left(\frac{x_{1}^{1/N}-x_{2}^{1/N}}{z_{0}}\right)
(x_{1}-x_{2})^{k}
\bar{Y}_{W}(u,x_{1})\bar{Y}_{W}(v,x_{2})\nonumber\\
& &\hspace{1cm}-z_{0}^{-1}\delta\left(\frac{x_{2}^{1/N}-x_{1}^{1/N}}{-z_{0}}\right)
(x_{1}-x_{2})^{k}\bar{Y}_{W}(v,x_{2})
\bar{Y}_{W}(u,x_{1})\nonumber\\
&=&x_{1}^{-1/N}\delta\left(\frac{x_{2}^{1/N}+z_{0}}{x_{1}^{1/N}}\right)
(x_{1}-x_{2})^{k}
Y_{W}(Y(\Delta(x_{1})u,z_{0})\Delta(x_{2})v,x_{2}^{1/N}).
\end{eqnarray}
It follows that there exists a nonnegative integer $k'\ge k$ such that
\begin{eqnarray}\label{elocality-twisted}
(x_{1}-x_{2})^{k'}[\bar{Y}_{W}(u,x_{1}),\bar{Y}_{W}(v,x_{2})]=0.
\end{eqnarray}
Next, we establish the weak associativity.
Let $w\in W$. Assume $u\in V^{j}$ for some $0\le j\le N-1$, 
i.e., $\sigma (u)=\omega_{N}^{j}u$.
{}From (\ref{elimit-property}), we have
$$x^{-j/N}\bar{Y}_{W}(u,x)\in (\End W)[[x,x^{-1}]].$$
Let $l$ be a nonnegative integer such that
$$x^{l-j/N}\bar{Y}_{W}(u,x)w\in W[[x]].$$
Using the commutation relation (\ref{elocality-twisted})
we get
\begin{eqnarray}
x_{1}^{l-j/N}(x_{1}-x_{2})^{k'}
\bar{Y}_{W}(u, x_{1})\bar{Y}_{W}(v, x_{2})w
\in W[[x_{1},x_{2}^{1/N}]][x_{2}^{-1/N}].
\end{eqnarray}
With (\ref{elocality-twisted}), 
{}from (\ref{ejacobi-deform}) we get
\begin{eqnarray*}
& &x_{1}^{-1/N}\delta\left(\frac{x_{2}^{1/N}+z_{0}}{x_{1}^{1/N}}\right)
(x_{1}-x_{2})^{k'}
Y_{W}(Y(\Delta(x_{1})u,z_{0})\Delta(x_{2})v,x_{2}^{1/N})w\\
&=&\left(z_{0}^{-1}\delta\left(\frac{x_{1}^{1/N}-x_{2}^{1/N}}{z_{0}}\right)
-z_{0}^{-1}\delta\left(\frac{x_{2}^{1/N}-x_{1}^{1/N}}{-z_{0}}\right)
\right)\\
& &\hspace{2cm}\cdot \left((x_{1}-x_{2})^{k'}
\bar{Y}_{W}(u,x_{1})\bar{Y}_{W}(v,x_{2})w\right)\nonumber\\
&=&x_{1}^{-1/N}\delta\left(\frac{x_{2}^{1/N}+z_{0}}{x_{1}^{1/N}}\right)
\left((x_{1}-x_{2})^{k'}
\bar{Y}_{W}(u,x_{1})\bar{Y}_{W}(v,x_{2})w\right).
\end{eqnarray*}
In view of Remark \ref{rsubstitutions1},
substituting $z_{0}=(x_{2}+x_{0})^{1/N}-x_{2}^{1/N},x_{2}>>x_{0}$ we get
\begin{eqnarray}
& &x_{1}^{-1/N}\delta\left(\frac{(x_{2}+x_{0})^{1/N}}{x_{1}^{1/N}}\right)
\nonumber\\
& &\ \ \ \ \cdot (x_{1}-x_{2})^{k'}
Y_{W}(Y(\Delta(x_{1})u,(x_{2}+x_{0})^{1/N}-x_{2}^{1/N})
\Delta(x_{2})v,x_{2}^{1/N})w\nonumber\\
&=&x_{1}^{-1/N}\delta\left(\frac{(x_{2}+x_{0})^{1/N}}{x_{1}^{1/N}}\right)
\left((x_{1}-x_{2})^{k'}
\bar{Y}_{W}(u,x_{1})\bar{Y}_{W}(v,x_{2})w\right).
\end{eqnarray}
For the expression on the left-hand side, using Proposition \ref{pbdm}, we have
\begin{eqnarray*}
& &x_{1}^{-1/N}\delta\left(\frac{(x_{2}+x_{0})^{1/N}}{x_{1}^{1/N}}\right)
x_{1}^{l-j/N}(x_{1}-x_{2})^{k'}\\
& &\ \ \ \ \cdot Y_{W}(Y(\Delta(x_{1})u,(x_{2}+x_{0})^{1/N}-x_{2}^{1/N})
\Delta(x_{2})v,x_{2}^{1/N})w\\
&=&x_{1}^{-1/N}\delta\left(\frac{(x_{2}+x_{0})^{1/N}}{x_{1}^{1/N}}\right)
x_{0}^{k'}(x_{2}+x_{0})^{Nl-j}\\
& &\ \ \ \ \cdot Y_{W}(Y(\Delta(x_{2}+x_{0})u,(x_{2}+x_{0})^{1/N}-x_{2}^{1/N})
\Delta(x_{2})v,x_{2}^{1/N})w\\
&=&x_{1}^{-1/N}\delta\left(\frac{(x_{2}+x_{0})^{1/N}}{x_{1}^{1/N}}\right)
x_{0}^{k'}(x_{2}+x_{0})^{Nl-j}
Y_{W}(\Delta(x_{2})Y(u,x_{0})v,x_{2}^{1/N})w\\
&=&x_{1}^{-1/N}\delta\left(\frac{(x_{2}+x_{0})^{1/N}}{x_{1}^{1/N}}\right)
x_{0}^{k'}(x_{2}+x_{0})^{Nl-j}\bar{Y}_{W}(Y(u,x_{0})v,x_{2})w.
\end{eqnarray*}
Then 
\begin{eqnarray*}
& &x_{1}^{-1/N}\delta\left(\frac{(x_{2}+x_{0})^{1/N}}{x_{1}^{1/N}}\right)
x_{0}^{k'}(x_{2}+x_{0})^{Nl-j}\bar{Y}_{W}(Y(u,x_{0})v,x_{2})w\\
&=&x_{1}^{-1/N}\delta\left(\frac{(x_{2}+x_{0})^{1/N}}{x_{1}^{1/N}}\right)
\left(x_{1}^{l-j/N}(x_{1}-x_{2})^{k'}
\bar{Y}_{W}(u,x_{1})\bar{Y}_{W}(v,x_{2})w\right)\\
& &\hspace{9cm}|_{x_{1}=x_{2}+x_{0}}\\
&=&x_{1}^{-1/N}\delta\left(\frac{(x_{2}+x_{0})^{1/N}}{x_{1}^{1/N}}\right)
\left(x_{1}^{l-j/N}(x_{1}-x_{2})^{k'}
\bar{Y}_{W}(u, x_{1})\bar{Y}_{W}(v, x_{2})w\right)\\
& &\hspace{9cm}|_{x_{1}=x_{0}+x_{2}}\\
&=&x_{1}^{-1/N}\delta\left(\frac{(x_{2}+x_{0})^{1/N}}{x_{1}^{1/N}}\right)
\left((x_{0}+x_{2})^{l-j/N}x_{0}^{k'}
\bar{Y}_{W}(u, x_{0}+x_{2})\bar{Y}_{W}(v, x_{2})w\right).
\end{eqnarray*}
Consequently, we have
\begin{eqnarray*}
& &(x_{0}+x_{2})^{l-j/N}x_{0}^{k'}
\bar{Y}_{W}(u, x_{0}+x_{2})
\bar{Y}_{W}(v,x_{2})w\\
&=&(x_{2}+x_{0})^{l-j/N}x_{0}^{k'}
\bar{Y}_{W}(Y(u, x_{0})v,x_{2})w,
\end{eqnarray*}
proving the weak associativity.
Therefore, $(W,\bar{Y}_{W})$ is a weak $\sigma$-twisted $V$-module.
\end{proof}

\end{document}